\documentclass[]{amsart}

\newtheorem{theorem}{Theorem}[section]

\theoremstyle{definition}

\newtheorem{example}[theorem]{Example}

\theoremstyle{remark}
\newtheorem{remark}[theorem]{Remark}

\numberwithin{equation}{section}

 \newcommand{\real}{{\bf R}}

\newcommand{\cc}{{\bf C}}

\newcommand{\bbar}{\left[ \begin{array}}
\newcommand{\ebar}{\end{array} \right] }
\newcommand{\bdm}{\begin{displaymath}}
\newcommand{\edm}{\end{displaymath}}
\newcommand{\beq}{\begin{equation}}
\newcommand{\beqa}{\begin{eqnarray}}
\newcommand{\beqas}{\begin{eqnarray*}}
\newcommand{\eeq}{\end{equation}}
\newcommand{\eeqa}{\end{eqnarray}}
\newcommand{\eeqas}{\end{eqnarray*}}
\newcommand{\dd}{\textup{d}}

\begin{document}

\title[Constant curvature submanifolds of
pseudo-{E}uclidean space]
{A loop group formulation for constant curvature submanifolds of
pseudo-{E}uclidean space}


\author{David Brander}
\address{Department of Mathematics\\ Faculty of Science\\Kobe University\\1-1, Rokkodai, Nada-ku, Kobe 657-8501\\ Japan}
\email{brander@math.kobe-u.ac.jp}

\author{Wayne Rossman}
\address{Department of Mathematics\\ Faculty of Science\\Kobe University\\1-1, Rokkodai, Nada-ku, Kobe 657-8501\\ Japan}
\email{wayne@math.kobe-u.ac.jp}


\begin{abstract}
We give a loop group formulation for the problem of isometric immersions
with flat normal bundle of
a simply connected pseudo-Riemannian manifold $M_{c,r}^m$, of dimension $m$,
constant sectional curvature $c \neq 0$, and signature $r$, 
into the pseudo-Euclidean
space $\real_s^{m+k}$, of signature $s\geq r$.  In fact
 these immersions are obtained canonically 
from the loop group maps corresponding
to isometric  immersions of the same manifold into
a pseudo-Riemannian sphere or hyperbolic space $S_s^{m+k}$ or
$H_s^{m+k}$, which have previously been studied.  A simple 
 formula is given for obtaining these immersions  from those loop group maps.
\end{abstract}

\maketitle

\section{Introduction} \label{intro}
Many special submanifolds can be formulated as maps into loop groups
which admit various techniques to produce or analyse solutions (see,
for example, \cite{bdpt2002}, \cite{terng2002}, \cite{fordywood},
 and associated references). Concerning the present article,
it was shown by Ferus and Pedit \cite{feruspedit1996} that isometric
 immersions 
with flat normal bundle, $Q_c^m \to Q_{\tilde{c}}^{m+k}$, between simply connected
Riemannian space forms,
where $c \neq \tilde{c}$ and $c \neq 0 \neq \tilde{c}$, admit a loop group
formulation. They come in natural families parameterised by a spectral parameter
$\lambda$ in either $\real^*$, $i\real^*$ or $S^1$, and the constant curvature
of the immersion varies with $\lambda$. They also showed how to produce 
many local solutions using the AKS theory, when $k \geq m-1$. If
$c < \tilde{c}$, then $k = m-1$ is the minimal codimension for even
a local isometric immersion, and in this critical codimension the
normal bundle is automatically flat. 
 
  In \cite{branderdorf}, it was shown that these loop group maps
can also be constructed from flat immersions, using a generalised DPW method.
In \cite{brander2}, each map was shown, moreover, to correspond to different
immersions into different target spaces, depending on whether the map was 
evaluated for values of the spectral parameter $\lambda$ in $\real^*$, $i\real^*$ or $S^1$, 
showing that various global isometric immersion problems are equivalent.

The immersions obtained for $\lambda$ in $i\real^*$ and $S^1$ shrink to a point
as $\lambda$ approaches the degenerate values $\pm i$, and so it is of interest to
find another way to interpret the loop group map at these values, to complete the
above picture. On the other hand, the case $\tilde{c} = 0$ has still not
been given a loop group formulation, and therefore is also of interest.
In this note we will simultaneously solve both of these problems.

We first prove, in Theorem \ref{prop1}, that isometric immersions
with flat normal bundle from a simply connected pseudo-Riemannian manifold
$M_{c,r}^m$, of constant curvature $c \neq 0$ and signature $r$,
 to a pseudo-Riemannian sphere or hyperbolic space, 
$S_s^{m+k}$ or $H_s^{m+k}$, of signature $s$, correspond, in a very natural way,
to isometric immersions with flat normal bundle
of the same manifold, $M_{c,r}^m$, into a pseudo-Euclidean space $\real_s^{m+k}$.
We give a simple, coordinate free, proof of this  result, which generalises 
to an arbitrary situation what had, in effect, been shown earlier \cite{tenetal} for the special 
case that the space being immersed is Riemannian, the codimension is $m-1$, and  such that the principle normal curvatures are never equal to $c - \tilde{c}$, which guaranteed the existence
 of principle coordinates.   
 
In Sections \ref{loopsec1} and \ref{loopsec}, we look at the loop group formulation mentioned
 above, which generalises easily to pseudo-Riemannian space forms of
arbitrary signatures.
Previously,  Cie\'{s}li\'{n}ski and
Aminov \cite{amci}, had shown locally and 
for the special case of isometrically immersing
the hyperbolic space $H^m$ into the sphere $S^{2m-1}$, 
that, by allowing the target space to grow so
that the curvature induced on the immersion remains fixed as the parameter $\lambda$ varies, one obtains, in the limit as the radius of the spherical
target space approaches infinity, an isometric immersion of $H^m$ in
Euclidean space $E^{2n-1}$.  They also gave a formula of Sym type for the 
immersion into $E^{2n-1}$.

We prove,  Theorem \ref{theorem2},
that this result holds globally, and for arbitrary signatures and
codimension. The limit as the target space approaches pseudo-Euclidean
space corresponds to the above-mentioned degenerate spectral parameter
values $\lambda = \pm i$, and the immersion into pseudo-Euclidean 
space obtained at $\lambda = i$ is the same  as the one given in
Theorem \ref{prop1}. We also prove that, conversely, every 
isometric immersion with flat normal bundle 
$M_{c,r}^{m} \to \real_{s}^{m+k}$ is associated 
to one of these loop group maps.

The Sym-type formula (\ref{sym}) allows one 
to obtain the immersion into pseudo-Euclidean
space directly from the loop group map.   Thus the loop group methods
for producing constant curvature immersions into pseudo-Riemannian
spheres and hyperbolic spaces, which have already been developed in 
 \cite{feruspedit1996} and \cite{branderdorf},
 automatically and explicitly produce isometric immersions
of the same manifold into pseudo-Euclidean space, with the same codimension.

\section{Constant curvature immersions with flat normal bundle into
 pseudo-Riemannian space forms}  \label{corsection}

Let $\real^n_s$ denote the pseudo-Euclidean space $\real^n$ with metric
of signature $s$, and $S^n_s$ denote the unit sphere in $\real^{n+1}_s$.
The isometry group of $S^n_s$ is $SO_s(n+1,\real) = \{ A \in GL(n+1,\real) ~|~ A^t J A = J \}$,
where $J$ is a diagonal matrix,  $\textup{diag}(\epsilon_1,...,\epsilon_{n+1})$, whose entries
$\epsilon_i$ are all $\pm1$,
and where $s$ of these entries are negative.
Let $ASO_s(n,\real)$ denote the group of pseudo-Euclidean motions of $\real^n_s$,
 that is the subgroup of $GL(n+1,\real)$ consisting of
matrices of the form 
\bdm
\bbar  {cc} T & a\\ 0 & 1 \ebar,
\edm
where $T \in SO_s(n,\real)$ and $a$ is a column vector.

Let $M$ be a  manifold of dimension $m$, and $\hat{f}: M \to \real_s^{m+k}$
an immersion such that the pull-back metric has signature $r$.
 An \emph{adapted
frame} for $\hat{f}$ is a map, $\hat{F}: \to ASO_s(m+k)$,
which has the form 
\bdm
\hat{F} = \bbar {cc} T & \hat{f} \\ 0 & 1 \ebar,
\edm
where $T = [\hat{e}_1,...,\hat{e}_m,\hat{n}_1,...,\hat{n}_k]$, and the column vectors $\hat{e}_i$
and $\hat{n}_j$ span the tangent and normal spaces respectively to the image of
$\hat{f}$.  We fix the matrix $J$ 
which defines $SO_s(n,\real)$ to be of the form
\bdm
J = \bbar {cc} J_1 & 0  \\ 0 & J_2 \ebar,
\edm
where $J_1$ is $m \times m$ and of signature $r$, and $J_2$ is $k \times k$ of signature
$s-r$.  $J_1$ encodes the signature of the induced metric.

 The Maurer-Cartan form for $\hat{F}$ is the pull-back of the Maurer-Cartan
form on $ASO_{s}(m+k)$, namely,
\beq \label{mcf}
\hat{A} = \hat{F}^{-1} \dd \hat{F} = \bbar {ccc}  \hat{\omega} &    \hat{\beta} & \hat{\theta} \\
                   \hat{\alpha} &  \hat{\eta} & 0 \\
                0 & 0 & 0 \ebar,
\eeq
where, if $\Omega (M)$ denotes the vector space of real valued 1-forms on
$M$, then $\hat{\omega} = [\hat{\omega}^i_j] \in so_r(m) \otimes \Omega (M)$
 and  $\hat{\eta} = [\hat{\eta}^i_j] \in so_{s-r}(k) \otimes \Omega (M)$ 
are the connections for the
tangent and normal bundles respectively,
$\hat{\beta} = [\hat{\beta}^i_j]$ is the second fundamental form,
$\hat{\alpha} = -J_2\hat{\beta}^t J_1$,
 and the components of
$\hat{\theta} = [\hat{\theta}^1,...,\hat{\theta}^m]^t$ make up the coframe dual
to $\hat{e}_i$.  To verify this, one uses the fact that 
${T}^{-1} = J {T}^t J$ and checks that the following 
equations are satisfied:
\beqas
\dd \hat{f} = \sum_j \hat{\theta}^j \hat{e}_j,\\
\dd \hat{e}_i = \sum_j \hat{\omega}^j_i \hat{e}_j + \sum _j \hat{\alpha}^j_i \hat{n}_j,\\
\dd \hat{n}_i = \sum_j \hat{\beta}^j_i \hat{e}_j + \sum _j \hat{\eta}^j_i \hat{n}_j.\\
\eeqas
Given such an immersion $\hat{f}$, an adapted frame always exists locally.
If $X$ and $Y$ are  matrix-valued 1-forms, their wedge product is defined
to have components $(X\wedge Y)^i_j := \sum_k X^i_k \wedge Y^k_j$.
The Maurer-Cartan form satisfies the Maurer-Cartan equation
\beq \label{mce}
\dd \hat{A} + \hat{A} \wedge \hat{A} = 0.
\eeq
Conversely, given a 1-form of the form (\ref{mcf}), defined on a 
simply connected subset
$U \subset M$, such that the components $\hat{\theta}^i$ are
all linearly independent 1-forms, and which satisfies  (\ref{mce}), then $\hat{A}$ integrates
to an adapted frame $\hat{F}$ for an immersion $\hat{f}: U \to \real^{m+k}$.

Let $Q_s^{n}(\epsilon)$ denote the pseudo-Riemannian sphere $S_s^{n}$ 
and hyperbolic space  $H_s^{n}$, 
for $\epsilon = 1$
and $\epsilon = -1$ respectively. The pseudo-Riemannian hyperbolic
space is defined as $H_s^{n} := \{ x \in \real_{s+1}^{n+1} ~|~ x^t J x = -1\}$,
where $J$ is the metric on $\real_{s+1}^{n+1}$.
Now suppose we have an immersion $f : M \to Q_s^{m+k}(\epsilon)$,
 with induced metric of signature $r$.  Then one has the analogue of
the above description for an adapted frame $F: M \to SO_{s+\delta}(m+k+1)$, 
where $\delta = \frac{1}{2}(1-\epsilon)$, and this
time $J= \textup{diag}(J_1, J_2, \epsilon)$,
\bdm
{F} = [e_1,...,e_m, n_1,...,n_k,f],
\edm
and 
\beq \label{mcf2}
{A} = {F}^{-1} \dd {F} = \bbar {ccc}  \omega &  \beta &   \theta \\
               \alpha & \eta & 0 \\
               - \epsilon  \theta ^t J_1 & 0 & 0 \ebar.
\eeq           
Here $\omega$, $\beta$, $\theta$, $\eta$ and $\alpha$ all
have the same form and interpretation as their corresponding objects
in the affine case above.
 
 In any of  the above cases, constant sectional curvature $c$ and 
flatness of the normal bundle for the immersion are
characterised respectively by the following equations and their
analogues, replacing an object $x$ with $\hat{x}$  where appropriate:
\beqa
\dd \omega + \omega \wedge  \omega = c \theta \wedge \theta^t J_1,  
\label{cc} \\
\dd \eta + \eta \wedge \eta = 0. \label{fnb}
\eeqa

\begin{theorem}  \label{prop1}
Let $M$ be a simply connected manifold of dimension $m$.
\begin{enumerate}
\item
Let ${f}: M \to S_s^{m+k}$ be a smooth immersion with flat normal bundle, such that the induced
metric has signature $r$, and  constant curvature $c \in (-\infty,0) \cup (1,\infty)$.  Then 
there is, uniquely up to an action by the isometry group of $\real_s^{m+k}$,
 a canonically defined immersion  with flat normal bundle $\hat{f}: M \to \real_s^{m+k}$,
with the same induced metric. 
The same statement holds with the roles of $S_s^{m+k}$  and $\real_s^{m+k}$ reversed.
\item Statement (1) is also valid  substituting $H_s^{m+k}$ for 
 $S_s^{m+k}$
and  $(-\infty,-1) \cup (0,\infty)$ for $(-\infty,0) \cup (1,\infty)$.
\end{enumerate} 
\end{theorem}
\begin{proof}
Let $f: M \to Q_s^{m+k}(\epsilon)$ be the map from either the
first or the second case.
Fix a base point $p$ of $M$. Without loss of generality, we assume that ${f}(p) = [0,...,0,1]^t$. Choose an adapted frame ${F}$ for ${f}$ on a simply
connected neighbourhood $U$ of $p$, normalised to the identity at $p$. This frame is unique up
to right multiplication by a smooth map 
${G}: U \to SO_r(m, \real) \times SO_{s -r}(k,\real) \subset 
SO_{s+\delta}(m+k +1,\real)$, which has the form 
\beqa \label{g}
{G} = \textup{diag}(G_1,G_2,1), && G(p) = I.
\eeqa
  This corresponds to 
a change of orthonormal frames for the tangent and normal bundles, while  fixing the last 
column, ${f}$.

Let ${A}$ be the Maurer-Cartan form for ${F}$, with components labelled as in the
equation (\ref{mcf2}).  Now set 
\beq \label{ahat2}
\hat{A} = \bbar {ccc}  {\omega} & \frac{i\sqrt{\epsilon c}}{\sqrt{1-\epsilon c}} {\beta} &  {\theta} \\
     \frac{i\sqrt{\epsilon c}}{\sqrt{1-\epsilon c}}\alpha  &  {\eta} & 0 \\
                0 & 0 & 0 \ebar.
\eeq
The allowed ranges for $c$ ensure that the factor 
$\frac{i\sqrt{\epsilon c}}{\sqrt{1-\epsilon c}}$ is real, so $\hat{A}$
 is a 1-form with values  in the Lie algebra of $ASO_{s}(m+k,\real)$.
It is a straightforward computation to verify that ${A}$ satisfies the
integrability condition (\ref{mce}) together with the equations (\ref{cc}) and
(\ref{fnb}) if and only $\hat{A}$ does also.  In the computation, one uses
(\ref{cc}) to obtain the equivalence of the first diagonal components
of the Maurer-Cartan equations, which for $A$ and $\hat{A}$ are, respectively,
\beqas
\dd \omega + \omega \wedge  \omega + \beta \wedge  \alpha - \theta 
    \wedge \epsilon \theta^t J_1 = 0, \\
\dd \omega + \omega \wedge  \omega - \frac{\epsilon c}{1-\epsilon c}\beta \wedge  \alpha = 0,
\eeqas
and the equation (\ref{fnb}) is needed for the second diagonal 
components.  Thus constant curvature and flatness of the normal bundle
are essential here.

  Now we can integrate $\hat{A}$ on $U$ to
get a unique adapted frame $\hat{F}$ for the desired immersion $\hat{f}$, with the initial condition
$\hat{F}(p)= I$. The freedom for the choice of adapted frame for $\hat{F}$ is also post-multiplication
by a smooth map $\hat{G}: U \to SO_r(m, \real) \times SO_{s-r}(k,\real) \subset 
ASO_{s}(m+k,\real)$, which has the same form as (\ref{g}), and has exactly the
same effect on the Maurer-Cartan form $\hat{A}$, whether it is applied
\emph{first} to ${F}$, and then constructing $\hat{A}$ as prescribed above,
or whether it is applied to $\hat{F}$ \emph{after} the construction from
$F$ and then differentiating $\hat{F}$. 
 Hence the map $\hat{f}$ is uniquely determined
by our choice of normalisation point $p$, which corresponds to an action of 
the isometry group $ASO_s(m+k,\real)$.

For the global picture, one observes that for any point $q$ in $M$, there is a
simply connected neighbourhood $U_q$ of $q$, which contains $p$, and an
adapted frame ${F}_q$ on $U_q$, normalised at $p$.   Thus the same 
procedure can be carried out on $U_q$.  On the overlap, $U_q \cap U$, ${F}$ and
${F}_q$ differ only by right multiplication by a matrix of the form (\ref{g}),
which has already been taken into account in our construction of $\hat{f}$ described above.

The induced metric for both $f$ and $\hat{f}$
 is given in terms of the local coframe, 
$J_1 \theta =: [\theta_1,...,\theta_m]^t$,
by the formula 
$\theta^2_1 + ... + \theta_m^2$, and thus is identical.  Clearly the same argument holds
with the roles of the target spaces reversed.
\end{proof}

\section{The loop group formulation}  \label{loopsec1}
The loop group  formulation for isometric immersions of space forms given in
\cite{feruspedit1996} works also for the pseudo-Riemannian case. Here is
a brief outline of the formulation of Ferus and Pedit, the only difference
here being that we allow non-Riemannian signatures. The computations are
easily verified to go through exactly as in \cite{feruspedit1996}.

Let $M$ be a simply
connected pseudo-Riemannian space form, with constant curvature $c \neq 0$ and
of signature $r$, and fix a base point $p$ of $M$.
Given an isometric immersion with flat normal bundle, $f$, of $M$
 into the pseudo-Riemannian
sphere or hyperbolic space, $Q^{m+k}_s(\epsilon)$, and an adapted frame $F$, one inserts a complex parameter $\lambda$ 
into the Maurer-Cartan form (\ref{mcf2}), to obtain a family of 1-forms,
\beq \label{alambda}
 A_\lambda = \bbar {ccc}  {\omega} & \frac{\sqrt{\epsilon c}}{2\sqrt{1-\epsilon c}}(\lambda - \lambda^{-1}){\beta} &   \frac{\sqrt{\epsilon c}}{2}(\lambda + \lambda^{-1}){\theta} \\
               \frac{\sqrt{\epsilon c}}{2\sqrt{1-\epsilon c}}(\lambda - \lambda^{-1}) \alpha  & {\eta} & 0 \\
       - \epsilon \frac{\sqrt{\epsilon c}}{2}(\lambda + \lambda^{-1}){\theta} ^t J_1 & 0 & 0 \ebar.
\eeq       
The original Maurer-Cartan form is obtained at $\lambda_0 = \frac{1}{\sqrt{\epsilon c}}(1 + \sqrt{1-\epsilon c})$.
The assumptions that $f$ has constant curvature and flat normal bundle are 
equivalent to the assumption that
$A_\lambda$ satisfies the Maurer-Cartan equation (\ref{mce}) for all $\lambda$ in the 
punctured plane $\cc^*$. 
 Depending on the original curvature value $c$,
$A_\lambda$ is real for $\lambda$ in one of $i\real^*$, $\real^*$, or $S^1$, and
integrates to an adapted frame 
$F_\lambda = [e_1^\lambda,...,e_m^\lambda,n_1^\lambda,...,n_k^\lambda,f_\lambda]$ for a family of immersions with flat 
normal bundle, with constant curvature in one of the corresponding ranges
$(-\infty,0)$, $(0,1)$, or $(1,\infty)$, for the case $\epsilon =1$,
and their reflections about 0 for the case $\epsilon = -1$.
The family $F_\lambda$ is unique with the normalisation $F_\lambda(p) = I$.
 As mentioned in the introduction,
there are several methods for producing the loop group maps $F_\lambda$,
for $k \geq m-1$, using
techniques from integrable systems.

\begin{remark} \label{globalremark}
An important point for our discussion in the next section is that, 
even though a single global adapted frame may not exist for $f$,
one can nevertheless show, \cite{brander2}, that the family $f_\lambda$ 
is well defined globally on $M$. 
\end{remark}

The coframe of the immersion $f_\lambda$, is given, from (\ref{alambda}),
by 
\beq \label{coframe}
\theta_\lambda = \frac{\sqrt{\epsilon c}}{2}(\lambda + \lambda^{-1}){\theta},
\eeq
and the curvature tensor turns out to be given by the expression
\beqas
\dd \omega + \omega \wedge \omega &=& c \theta \wedge \theta ^t J_1\\
 &=& \frac{4 \epsilon}{(\lambda + \lambda^{-1})^2} \theta_\lambda \wedge \theta_\lambda ^t J_1.
\eeqas
Thus, the constant curvature, which varies
with $\lambda$, is given by the formula
\bdm
c_\lambda = \frac{4 \epsilon}{(\lambda + \lambda^{-1})^2}.
\edm
Because the original map, $f = f_{\lambda_0}$, is
 an immersion, it follows from (\ref{coframe}) that $f_\lambda$ is an 
immersion for all $\lambda \neq \pm i$.
However,  $f_{\pm i}$ maps $M$ to a single point, because its coframe is zero.
 If we take the normalisation $F_\lambda(p) = I$, for some $p \in M$, then 
\bdm
f_{\pm i} = [0,...,0,1]^t.
\edm

\section{Interpretation of the loop group map at $\lambda = \pm i$} \label{loopsec}
We seek an  interpretation of the map $F_\lambda$ at $\lambda = \pm i$.
If we scale the space that we are immersing into,
so that the curvature of the immersion is unchanged, instead of varying with $\lambda$,
then, as $\lambda$ approaches $\pm i$, we must have  the target space
approaching flat space, and so, in the limit, we hope to have a 
constant curvature immersion $\hat{f} : M \to \real_s^{m+k}$.  

To carry this out, set
\beq \label{ftilde}
 \tilde{f}_\lambda := \frac{2}{\sqrt{\epsilon c}(\lambda + \lambda^{-1})} f_\lambda.
\eeq
As with the expression (\ref{alambda}), it is easy to verify that $\tilde{f}_\lambda$ is real for values of $\lambda$ in the appropriate ranges
corresponding to $c$ and $\epsilon$, described above. Thus 
$\tilde{f}_ \lambda$ is a map from  $M$ into 
$Q_s^{m+k}({\epsilon |R|}) := \{ x \in \real_{s+\delta}^{m+k+1} ~|~ 
 x^t J x = \epsilon R^2\}$, where 
$R := \frac{2}{\sqrt{\epsilon c} (\lambda + \lambda ^{-1})}$. 
Now the image of $\tilde{f}$ can be identified with the image of
$f$ if we scale the ambient space $\real_{s+\delta}^{m+k+1}$ by a factor
of $R$. This has the effect of scaling the metric by a constant
conformal factor of ${R^2}$, and the curvature by $\frac{1}{R^2}$.
 Thus  $\tilde{f}_\lambda$ has constant curvature
\bdm
c_\lambda \frac{\epsilon c(\lambda + \lambda^{-1})^2}{4} = c.
\edm
Because all we have done is scale the target space, an adapted frame, 
$\tilde{F}_\lambda \in ASO_{s+\delta}(m+k+1,\real)$, for $\tilde{f}_\lambda : M \to \real_{s+\delta}^{m+k+1}$, is
\bdm
\tilde{F}_\lambda = \bbar {cccccccc} e_1^\lambda & ... & e_m^\lambda &  n_1^\lambda & ... & n_k^\lambda & f_\lambda & \tilde{f}_\lambda \\
       &0 & &  & 0 & &  0& 1 \ebar,
\edm
where $e_i^\lambda $ and $n_i^\lambda$ are the same as in the
 unscaled frame $F_\lambda$. 
The unscaled $f_\lambda$ is now the $(k+1)$'st unit normal vector. 
The Maurer-Cartan form for $\tilde{F}_\lambda$ is
\beq \label{atilde}
 \tilde{A}_\lambda = \bbar {cccc}  {\omega} & \frac{\sqrt{\epsilon c}}{2\sqrt{1-\epsilon c}}(\lambda - \lambda^{-1}){\beta} &   \frac{\sqrt{\epsilon c}}{2}(\lambda + \lambda^{-1}){\theta}&  {\theta} \\
                \frac{\sqrt{\epsilon c}}{2\sqrt{1-\epsilon c}}(\lambda - \lambda^{-1}) \alpha & {\eta} & 0 & 0\\
     - \epsilon \frac{\sqrt{\epsilon c}}{2}(\lambda + \lambda^{-1}){\theta} ^t J_1 & 0 & 0 & 0\\
               0 & 0 & 0 & 0\ebar.
\eeq  
The frame $\tilde{F}_\lambda$ is obtained by 
integrating this with the initial condition 
\bdm
\tilde{F}(p) = \bbar {cc} 
   I_{m+k+1} & [0,...,0,\frac{2}{\sqrt{\epsilon c}(\lambda + \lambda^{-1})}]^t \\
   0 & 1  \ebar.
\edm

Now $\tilde{F}_\lambda$ is not defined at $\lambda = \pm i$,
because $\tilde{f}_\lambda$ is not. However, $\tilde{A}_\lambda$ \emph{is}
 defined at $\lambda = \pm i$, and at that point reduces to 
\bdm
\tilde{A}_{\pm i} = \bbar {cccc} {\omega} &  \frac{(\pm i)\sqrt{\epsilon c}}{\sqrt{1-\epsilon c}}{\beta} & 0&  {\theta} \\
     \frac{(\pm i)\sqrt{\epsilon c}}{\sqrt{1-\epsilon c}} \alpha & {\eta} & 0 & 0\\
     0 & 0 & 0 & 0\\
               0 & 0 & 0 & 0\ebar.
\edm
If we apply a translation to our initial condition for
 $\tilde{f}_\lambda$, and integrate
$\tilde{A}_{\lambda}$ with the initial condition $\hat{F}(p) = I_{m+k+2}$,
we get an adapted frame $\hat{F}_\lambda$ which satisfies
\bdm
\hat{F}_{\pm i} =  \bbar {cccccccc} e_1^{\pm i} & ... & e_m^{\pm i} &  n_1^{\pm i} & ... & n_k^{\pm i} & [0,...,0,1]^t & \hat{f}_{\pm i} \\
       & 0 & & & 0 &&  0& 1 \ebar,
\edm
and we see that $\hat{f}_{\pm i}$ maps into the hyperplane perpendicular to
the vector $[0,...,0,1]^t$.  In other words, it is an immersion into $\real_s^{n+k} \subset \real_{s+\delta}^{n+k+1}$, with flat normal bundle and constant curvature $c$.  Comparing the 
Maurer-Cartan form of $\hat{F}_{\pm i}$ with (\ref{ahat2}), we have shown that 
 the interpretation of the loop group map $F_\lambda$ at $i$ is just the 
 immersion into pseudo-Euclidean space obtained from Theorem \ref{prop1}.

Finally, since $\tilde{f}_\lambda = \hat{f}_\lambda +
   [0,...0,\frac{2}{\sqrt{\epsilon c}(\lambda + \lambda^{-1})}]^t$, we obtain, using
(\ref{ftilde}), the formula
\beq \label{fhatl}
\hat{f}_\lambda = \frac{2}{\sqrt{\epsilon c}(\lambda + \lambda^{-1})} (f_\lambda - [0,...,0,1]^t).
\eeq
Setting 
 $\mu := \lambda + \lambda^{-1}$ and $g(\mu) := f_{\lambda(\mu)}$,
 we obtain a Sym-type formula for  $\hat{f}_{\pm i}$:
\beqas
\hat{f}_{\pm i} &=& \frac{2}{\sqrt{\epsilon c}}
  \lim_{\mu \to 0} \frac{g(\mu)  - g(0)}{\mu}\\
&=&  \frac{2}{\sqrt{\epsilon c}}\frac{\partial \lambda}{\partial \mu}\frac{\partial}{\partial \lambda} f_\lambda \Big |_{\lambda = \pm i}\\
&=& \frac{1}{\sqrt{\epsilon c}} \frac{\partial}{\partial \lambda} f_\lambda\Big |_{\lambda = \pm i}.
\eeqas 
This formula is independent of the choice of adapted frame, $F$, and is
therefore,  by Remark \ref{globalremark}, valid globally.

Conversely, given an isometric immersion with flat normal bundle
 $f: M \to \real_{s}^{m+k}$, it follows from the converse part of
 Theorem \ref{prop1} that there are unique loop group maps $F_\lambda$,
normalised at $p$, with Maurer-Cartan forms of the form (\ref{alambda}),
corresponding to $f$.  We summarise this discussion as:
\begin{theorem} \label{theorem2}
Let $f: M_{c,r}^m \to Q_s^{m+k}(\epsilon)$ be an 
isometric immersion
with flat normal bundle of a simply connected $m$-dimensional 
pseudo-Riemannian manifold $M_{c,r}^m$, 
of signature $r$, constant curvature $c$ and with base point $p$,
 into the pseudo-Riemannian sphere or hyperbolic space, for $\epsilon =1$,
or $\epsilon =-1$ respectively.  Suppose that 
$c\in (-\infty,0) \cup (1,\infty)$ if $\epsilon =1$, and
$c \in (-\infty,-1) \cup (0,\infty)$ if $\epsilon =-1$.

Let $f_\lambda$ be the associated family of immersions given by
the last column of the frame  obtained by
integrating the 1-form defined by (\ref{alambda}), normalised
at $p$. Then the corresponding isometric immersion 
$\hat{f}: M_{c,r}^m \to \real_s^{m+k}$ from Theorem \ref{prop1} is given by
the formula:
\beq  \label{sym}
\hat{f} = \frac{1}{\sqrt{\epsilon c}} \, \pi_{m+k}  \left\{ \frac{\partial}{\partial \lambda} f_\lambda\big |_{\lambda =  i} \right\},
\eeq
where $\pi_{m+k}$ is the projection onto the first $m+k$ coordinates.

Conversely, every isometric immersion with flat normal bundle $f: M_{c,r}^{m} \to \real_{s+\delta}^{m+k}$ is obtained in this way.
\end{theorem}
\begin{remark} The formula (\ref{fhatl}) for $\hat{f}_\lambda$
gives a continuous deformation
of the original immersion $f -[0,...,0,1]^t$ into (the 
displaced)
$Q_s^{m+k}(\epsilon)$, obtained at $\lambda = \lambda_0$,
through to the immersion into $\real_s^{m+k} \subset \real_s^{m+k+1}$,
obtained at $\lambda = \pm i$.
\end{remark}
\begin{example}
As a simple test case, here is an example of a family,
from the loop group construction described above, of immersions
 into $H^3_1$ of the de Sitter spaces $S_{c_\lambda ,1}^2$ 
with constant sectional curvature 
$c_\lambda \in (0, \infty)$,
 for  values of $\lambda$ in $i \real ^*\setminus \{\pm i\}:$
\beqas
f_\lambda (u,v) = [ -i a \cosh u \sin v,
          ~ - i a \sinh u,
       ~ ab(1- \cos v \cosh u),
        ~  a^2 \cos v  \cosh  u  - b^2 ]^t, \\
F_\lambda = \bbar {cccc} \cos v & \sin v \sinh u & -ib\cosh u \sin v & -ia \cosh u \sin v \\
   0 & \cosh u & -ib \sinh u & -i a \sinh u \\
  ib \sin v & -ib \cos v \sinh u & a^2 -b^2 \cos v \cosh u & ab (1- \cos v \cosh u) \\
-ia \sin v & ia \cos v \sinh u & ab(\cos v \cosh u -1)&  a^2 \cos v \cosh u -b^2 
  \ebar, \\
a := \frac{1}{2}(\lambda + \lambda^{-1}) \hspace{1cm} 
  b := \frac{1}{2}(\lambda - \lambda^{-1}).
\eeqas 
 For any $c \in (0,\infty)$, we apply the
above formula at $\lambda =i$ to obtain 
\bdm
\hat{f}(u,v) = \frac{1}{\sqrt{c}}[-\cosh u  \sin v,~ -\sinh  u, ~1 -\cos v \cosh u]^t,
\edm
an embedding of the de Sitter space $S_{c,1}^2$ of constant curvature $c$ into
 $\real_1 ^3$.

Analogous test cases with other signatures
can be similarly constructed. A Riemannian example can be computed
using the example in \cite{brander2}.
\end{example}

\bibliographystyle{amsplain} 

\bibliography{mybib}

\providecommand{\bysame}{\leavevmode\hbox to3em{\hrulefill}\thinspace}
\providecommand{\MR}{\relax\ifhmode\unskip\space\fi MR }
\providecommand{\MRhref}[2]{%
  \href{http://www.ams.org/mathscinet-getitem?mr=#1}{#2}
}
\providecommand{\href}[2]{#2}
\begin{thebibliography}{1}

\bibitem{tenetal}
JL~Barbosa, W~Ferreira, and K~Tenenblat, \emph{Submanifolds of constant
  sectional curvature in pseudo-{R}iemannian manifolds}, Ann. Global Anal.
  Geom. \textbf{14} (1996), 381--401.

\bibitem{brander2}
D~Brander, \emph{Curved flats, pluriharmonic maps and constant curvature
  immersions into pseudo-{R}iemannian space forms}, arxiv preprint:0610058
  (2006).

\bibitem{branderdorf}
D~Brander and J~Dorfmeister, \emph{The generalized {DPW} method and an
  application to isometric immersions of space forms}, arxiv preprint:0604247
  (2006).

\bibitem{bdpt2002}
M~Br\"{u}ck, X~Du, J~Park, and C~Terng, \emph{The submanifold geometries
  associated to {G}rassmannian systems}, Mem. Amer. Math. Soc. 155 (2002),
  no.~735, viii + 95 pp.

\bibitem{amci}
JL~Cie\'{s}li\'{n}ski and Y~A Aminov, \emph{A geometric interpretation of the
  spectral problem for the generalized sine-{G}ordon system}, J. Phys. A: Math.
  Gen. \textbf{34} (2001), L153--L159.

\bibitem{feruspedit1996}
D~Ferus and F~Pedit, \emph{Isometric immersions of space forms and soliton
  theory}, Math. Ann. \textbf{305} (1996), 329--342.

\bibitem{fordywood}
A~P Fordy and J~C Wood (eds.), \emph{Harmonic maps and integrable systems},
  Aspects of mathematics, Vieweg, 1994.

\bibitem{terng2002}
C~Terng, \emph{Geometries and symmetries of soliton equations and integrable
  elliptic equations}, arXiv preprint:0212372 (2002).

\end{thebibliography}

\end{document}